\documentclass[12pt,a4paper,reqno]{amsart}
\usepackage{amsmath,amsthm,amssymb}

\usepackage{color}
\usepackage[dvips]{graphicx}
\theoremstyle{plain}
\newtheorem{theorem}{Theorem}
\newtheorem{lemma}{Lemma}

\newtheorem{corollary}{Corollary}
\newtheorem{definition}{Definition}
\theoremstyle{definition}

\usepackage [right=2.5cm, left=2.5cm, top=2.5cm, bottom=2.5cm]{geometry}

\begin{document}
\title{On weakly $m$-convex sets, $m=1,\ldots, n-1$,  in $n$-dimensional Euclidean spaces}
\author{Tetiana M. Osipchuk}
\address{Institute of Mathematics of NASU, 01004 Ukraine, Kiev-4\\
3, Tereschenkivska st.}
\email{osipchuk.tania@gmail.com}
\begin{abstract}
We study general topological properties of open and closed weakly $m$-convex sets, $m=1,\ldots, n-1$, in the real $n$-dimensional Euclidean space $\mathbb{R}^n$, $n\geqslant2$. An open set in $\mathbb{R}^n$ is called \textit{weakly $m$-convex} if, through every boundary point of the set, there passes an $m$-dimensional plane that does not intersect the set. A closed set in $\mathbb{R}^n$ is \textit{weakly $m$-convex} if it can be approximated from the outside by a countably decreasing family of open weakly $m$-convex sets. The set of all points in the complement of a given subset of $\mathbb{R}^n$ through each of which every $m$-dimensional plane intersects the subset is called the \textit{$m$-nonconvexity-point set} corresponding to the given subset.  We prove that the nonempty $m$-nonconvexity-point set corresponding to an open weakly $m$-convex set is itself open
and weakly $m$-convex, and that the interior of every closed weakly $m$-convex set with nonempty interior is weakly $m$-convex for every $m=1,\ldots, n-1$.
\end{abstract}


\begingroup
\def\uppercasenonmath#1{}
\maketitle
\endgroup
\section{Introduction}
We study generalized convex sets in the real $n$-dimensional Euclidean space $\mathbb{R}^n$, $n\geqslant2$, known as weakly $m$-convex sets, where $m=1,2,\ldots,n-1$.  The notion was coined by Yu. B. Zelinskii and I. V. Momot in 2001 \cite{Zel01} and has since been further developed in a number of works.

An open subset $E\subset\mathbb{R}^n$ is called \textit{weakly $m$-convex}  if, through every boundary point of $E$, there passes an $m$-dimensional plane that does not intersect $E$. A closed subset of $\mathbb{R}^n$ is \textit{weakly $m$-convex} if it can be approximated from the outside by a countably decreasing family of open weakly $m$-convex sets.  For a subset $E\subset\mathbb{R}^n$ a point in $\mathbb{R}^n\setminus E$ is called an \textit{$m$-nonconvexity point} of $E$ if every $m$-dimensional plane passing through the point intersects $E$. The collection of all $m$-nonconvexity points of  $E\subset\mathbb{R}^n$ is called the \textit{$m$-nonconvexity-point set} corresponding to $E$ and is denoted by $E_m^{\triangle}$.

It is known that, for every open weakly $m$-convex subset $E\subset\mathbb{R}^n$, $n\geqslant2$, where $m=1, n-1$, with $E_m^{\triangle}\ne\varnothing$, the set $E_m^{\triangle}$ is open and weakly $m$-convex, \cite{Osi2022_2}, \cite{Osi2023}, \cite{Osi2025}, \cite{Osi2026_1}, \cite{Osi2026_3}. Furthermore, the nonempty interior of a closed weakly $m$-convex set in  $\mathbb{R}^n$, $n\geqslant2$, is weakly $m$-convex for $m=1, n-1$ \cite{Osi2023}, \cite{Osi2025}, \cite{Osi2026_2}, \cite{Osi2026_3}. In this paper, we establish these results for the remaining cases $n\geqslant4$ and $1<m<n-1$. The proofs also cover the endpoint cases $m=1$ and $m=n-1$, for which the corresponding results are already known.

\section{Definitions and Preliminaries}

An $m$-dimensional affine subspace of $\mathbb{R}^n$, $n\geqslant1$, where  $m=0,1,2,\ldots,n-1$, is called an {\bf\emph{$m$-dimensional plane}}. An $(n-1)$-dimensional plane in $\mathbb{R}^n$ is also called a {\bf\emph{hyperplane}}.

We use the following standard notations. For a subset $E\subset\mathbb{R}^n$, let $\overline{E}$
be its closure, $\mathrm{Int}\, E$ be its interior, and $\partial
E=\overline{E}\setminus\mathrm{Int}\, E$ be its boundary. 

\begin{definition} [\cite{Zel01}]\label{def3}
An open subset  $E\subset\mathbb{R}^n$, $n\geqslant2$, is called {\bf\emph{weakly $m$-convex}},  $m=1,2,\ldots,n-1$, if, for every point $x\in\partial E$, there exists
an $m$-dimensional plane $L$ such that $x\in L$ and $L\cap E=\varnothing$.
\end{definition}

\begin{definition}\label{def2}
A set $E$ is said to {\bf\emph{be approximated from the outside}} by a family of open sets $E_k$,
$k=1,2,\ldots$, if $\overline{E}_{k+1}$ is contained in $E_k$, and $E=\bigcap\limits_kE_k$.
\end{definition}

It can be proved that any set approximated from the outside by a
family of open sets is closed.

\begin{definition}[\cite{Zel01}]\label{def4}
A closed subset $E\subset\mathbb{R}^n$, $n\geqslant2$, is called {\bf\emph{weakly $m$-convex}},  $m=1,2,\ldots,n-1$, if it can be approximated
from the outside by a family of open weakly $m$-convex sets.
\end{definition}
Among the closed weakly
$m$-convex sets there are also those with empty interior:
$$
E=\overline{E}=\overline{E}\setminus \mathrm{Int}\, E=\partial E.
$$

The class of open and closed weakly $m$-convex sets in $\mathbb{R}^n$ is denoted by $\mathbf{WC^n_m}$.

It is not difficult to verify that any open convex subset $E\subset\mathbb{R}^n$, $n\geqslant2$, is weakly $m$-convex, $m=1,2,\ldots,n-1$. Indeed, through every boundary point of a convex set, there passes a supporting hyperplane for the set; see Theorem~3.2, \cite{Leh}. Since $E$ is open, the supporting hyperplane does not intersect $E$. Thus, $E$ is weakly $(n-1)$-convex and, therefore, weakly $m$-convex for every $1\leqslant m\leqslant n-1$.

Similarly,  any closed convex subset of $\mathbb{R}^n$ is weakly $m$-convex, $1\leqslant m\leqslant n-1$, since it can be approximated from the outside by a decreasing family of open convex neighborhoods.

However, there exist disconnected, and hence non-convex,  weakly $m$-convex sets for every $m=1,2,\ldots,n-1$.

\begin{definition}\label{def5}
For a subset $E\subset\mathbb{R}^n$, $n\geqslant2$, a point $x\in \mathbb{R}^n\setminus E$ is called an {\bf\emph{$m$-nonconvexity point}}, $m=1,2,\ldots,n-1$,  if every $m$-dimensional plane passing through $x$ intersects $E$. The set of all $m$-nonconvexity points of $E$ with respect to a fixed $m\in\{1,2,\ldots,n-1\}$ is called the {\bf\emph{$m$-nonconvexity-point set}} corresponding to $E$, and  is denoted by $E_m^{\triangle}$.
\end{definition}

\section{Main results}
Given two points $x,y\in\mathbb{R}^n$, we denote by $xy$ the
open line segment joining $x$ and $y$, by $\|x-y\|$ the length
of this segment, and by $\langle x,y\rangle$ the Euclidean
inner product of the vectors $x$ and $y$. We also set
$$
S^{n-1}:=\{x\in\mathbb{R}^n:\|x\|=1\}
$$
and
$$
U(y,r):=\{x\in\mathbb{R}^n:\|x-y\|<r\},
\quad y\in\mathbb{R}^n,\quad r>0.
$$
Let
$$
\mathcal L:=\{(a_1,\ldots,a_m)\in (\mathbb{R}^n)^m: \,\,\mbox{vectors}\,\,a_1,\ldots,a_m\in\mathbb{R}^n \,\,\mbox{are linearly independent}\}.
$$
For $y\in\mathbb{R}^n$ and
$(a_1,\ldots,a_m)\in \mathcal L$, where $1\leqslant m\leqslant n-1$,
we define
$$
L(y,a_1,\ldots,a_m):=
\left\{y+\sum_{j=1}^{m}\alpha_j a_j:
\alpha_j\in\mathbb{R},\ j=1,\ldots,m\right\}.
$$
Let $O$ denote the origin of $\mathbb{R}^n$, and write
$$
L(a_1,\ldots,a_m):=L(O,a_1,\ldots,a_m).
$$

For $1\leqslant m\leqslant n-1$, let $G(n,m)$ denote the Grassmannian of all $m$-dimensional linear subspaces of $\mathbb{R}^n$. We identify each subspace $L\in G(n,m)$ with the matrix $P_L$ of the orthogonal projection onto $L$. Thus, $G(n,m)$ is identified with the set
$$
\mathcal{P}_{n,m}:=\left\{P\in\mathbb{R}^{n\times n}: P^T=P,\,\, P^2=P,\,\, \operatorname{rank}P=m\right\}.
$$
We equip $G(n,m)$ with the topology induced by the standard topology of $\mathbb{R}^{n\times n}$ under this identification; see \cite{Ben}.
\begin{lemma}\label{lem1}
Let $U(p,r)$ be an open ball in $\mathbb{R}^n$, where $n\geqslant2$.
Then the set
$$
V:=\left\{L\in G(n,m):L\cap U(p,r)\ne\varnothing\right\}
$$
is open in $G(n,m)$.
\end{lemma}
\begin{proof}
Define the function $d_p:G(n,m)\to[0,+\infty)$ by
$$
d_p(L):=\inf_{x\in L}\|p-x\|.
$$
Then
$$
V=\{L\in G(n,m):d_p(L)<r\}.
$$
The function $d_p$ is continuous on $G(n,m)$. Indeed, by the
properties of orthogonal projections, we have
$$
d_p(L)=\|(I-P_L)p\|,
$$
where $I$ is the identity operator on $\mathbb R^n$ and $P_L$
is the orthogonal projection onto $L$. Consequently, $I-P_L$
is the orthogonal projection onto the orthogonal complement
$L^\perp$ of $L$.

By the definition of the topology on $G(n,m)$, the map
$L\mapsto P_L$ is continuous. Since the map
$P\mapsto (I-P)p$ and the Euclidean norm are also continuous,
it follows that $d_p$ is continuous on $G(n,m)$.

Finally,
$$
V=d_p^{-1}((-\infty,r)).
$$
Since $(-\infty,r)$ is open in $\mathbb R$ and $d_p$ is
continuous, $V$ is open in $G(n,m)$.
\end{proof}
\begin{lemma}\label{lem2}
Let $E\subset \mathbb{R}^n$, $n\geqslant2$, be an open subset such that $E_m^{\triangle}\ne\varnothing$, where $m\in\{1,2,\ldots,n-1\}$. Let $y\in E_m^{\triangle}$. Then, for every $L (y,a_1,\ldots,a_m)$, $(a_1,\ldots,a_m) \in \mathcal L$, there exists an open ball $U(x(y,a_1,\ldots,a_m),d(y))\subset E$, where $x(y,a_1,\ldots,a_m)\in L (y,a_1,\ldots,a_m)\cap E$ and $d(y)>0$ depends only on $y$, not on $a_1,\ldots,a_m$.
\end{lemma}
\begin{proof}
Fix an arbitrary point $y\in E_m^{\triangle}$. Without loss of generality, suppose that $y=O$.
Since $O\in E_m^{\triangle}$, every $m$-dimensional plane $L (a_1,\ldots,a_m)$, $(a_1,\ldots,a_m) \in \mathcal L$, intersects $E$. For each $L (a_1,\ldots,a_m)$, choose a point $p(a_1,\ldots,a_m)\in L (a_1,\ldots,a_m)\cap E$.
Since $E$ is open, there exist open balls $U_{a_1,\ldots,a_m}:=U\left(p(a_1,\ldots,a_m),r(a_1,\ldots,a_m) \right)$ such that $\overline{U_{a_1,\ldots,a_m}}\subset E$.

Let
$$V_{a_1,\ldots,a_m}:=\left\{ L \in G(n,m): L\cap U_{a_1,\ldots,a_m} \ne\varnothing\right\},\quad (a_1,\ldots,a_m) \in \mathcal L.$$

By Lemma~\ref{lem1}, the sets $V_{a_1,\ldots,a_m}$ are open in $G(n,m)$, moreover, $L(a_1,\ldots,a_m) \in V_{a_1,\ldots,a_m}$. Then
$$
\bigcup\limits_{(a_1,\ldots,a_m) \in \mathcal L} V_{a_1,\ldots,a_m}
$$
is a cover of $G(n,m)$.
Since $G(n,m)$ is compact; see \cite{Ben}, by the Heine-Borel theorem, there exists a subcover
$$
\bigcup\limits_{j\in M}
V_{a_1^j,\ldots,a_m^j},\quad (a_1^j,\ldots,a_m^j) \in \mathcal L,\quad j\in M,\quad M\,\,\mbox{is finite},
$$
of $G(n,m)$. Consequently, any $m$-dimensional plane $L(a_1,\ldots,a_m)$, $(a_1,\ldots,a_m)\in \mathcal L$, intersects $\bigcup\limits_{j\in M}\overline{U_{a_1^j,\ldots,a_m^j}}$.

Let $E_i$, $i\in N\subseteq \mathbb{N}$, be the components of $E$. Then for any $j\in M$ there exists $i(j)\in N$ such that $\overline{U_{a_1^j,\ldots,a_m^j}}\subset E_{i(j)}$.

Consider the distance functions
$$d_i(x):=\inf\limits_{x^0\in \partial E_i}\|x-x^0\|,\quad x\in E_i,\quad i\in N.$$
They are continuous in the domains $E_i$, $i\in N$. Then their restrictions to the compacts $\overline{U_{a_1^j,\ldots,a_m^j}}$ attain their minimum values $d_{j}>0$ on these compacts, i.e.,
$$
d_{j}:=\min\limits_{x\in \overline{U_{a_1^j,\ldots,a_m^j}}}d_{i(j)}(x),\quad j\in M.
$$
Since $M$ is finite, there exists
$$
d:=\min\limits_{j\in M}d_j>0.
$$

\begin{figure}[h]
\centering
\includegraphics[width=8 cm]{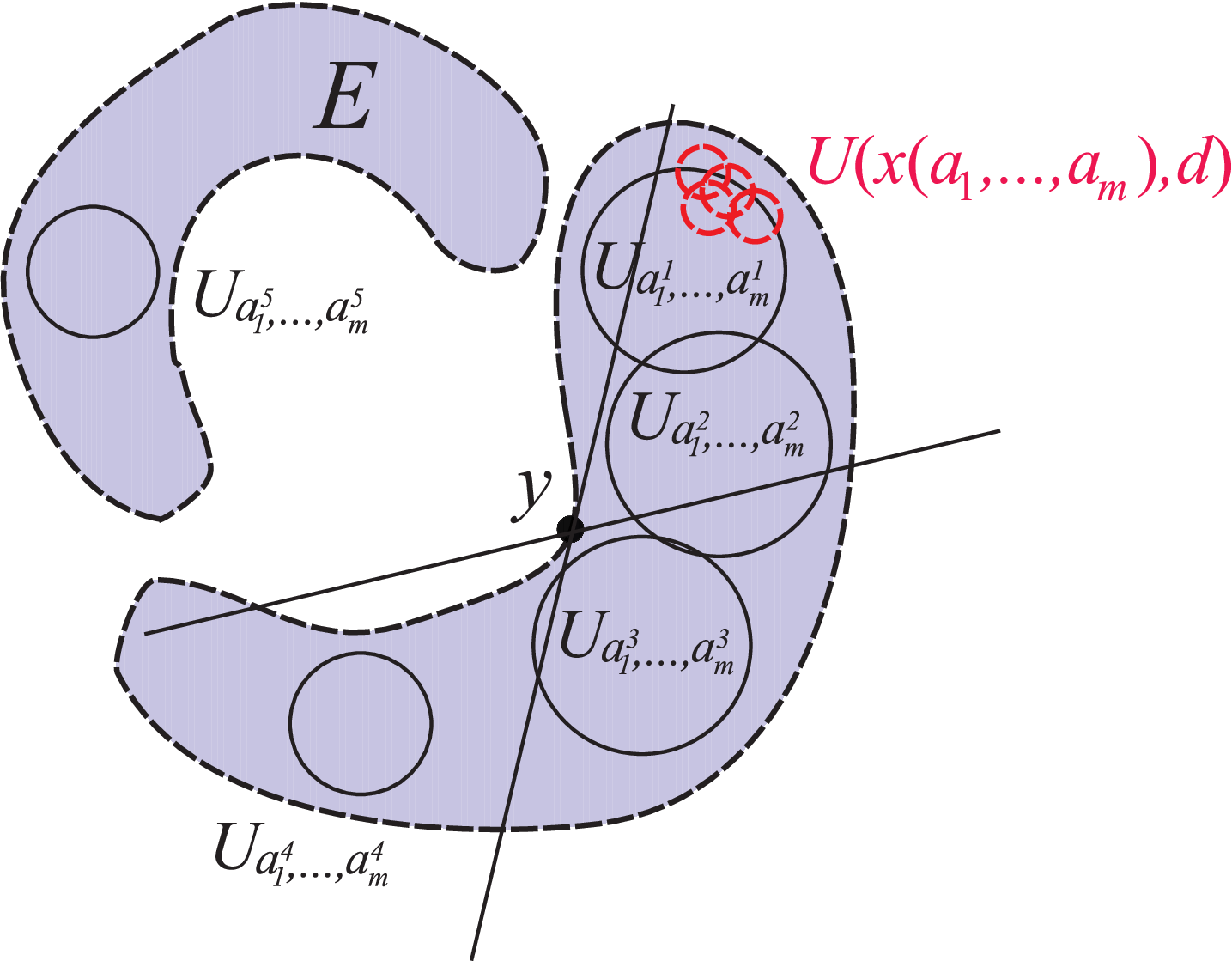}
\caption{}\label{Fig4}
\end{figure}

By the construction of the finite subcover, for every $(a_1,\ldots,a_m)\in\mathcal L$, there exists $j\in M$ such that $L(a_1,\ldots,a_m)\cap U_{a_1^j,\ldots,a_m^j}\ne\varnothing$.
Choose a point
$x(a_1,\ldots,a_m)\in L(a_1,\ldots,a_m)\cap U_{a_1^j,\ldots,a_m^j}$.
Then,
$$
U(x(a_1,\ldots,a_m),d)\subset E.
$$
Thus, the required assertion follows.
\end{proof}
\begin{theorem}
Let an open subset $E\subset \mathbb{R}^n$, $n\geqslant2$, belong to the class $\mathbf{WC^n_m}$, and $E_m^{\triangle}\ne\varnothing$, where $m=1,2,\ldots,n-1$. Then $E_m^{\triangle}$ is open.
\end{theorem}
\begin{proof} Fix an arbitrary point $y\in E_m^{\triangle}$ and show that it is an inner point of $E_m^{\triangle}$.

Since $E$ is weakly $m$-convex, it follows that $y\not\in \partial E$. Then there exists a number $\varepsilon_1>0$ such that $U(y,\varepsilon_1)\subset \mathbb{R}^n\setminus \overline{E}$.

By Lemma~\ref{lem2}, for the fixed $y$ and every $(a_1,\ldots,a_m) \in \mathcal L$ there exist a point $x(a_1,\ldots,a_m) \in L (y,a_1,\ldots,a_m)\cap E$ and a constant $d>0$ such that $U(x(a_1,\ldots,a_m),d)\subset E$.

\begin{figure}[h]
\centering
\includegraphics[width=16 cm]{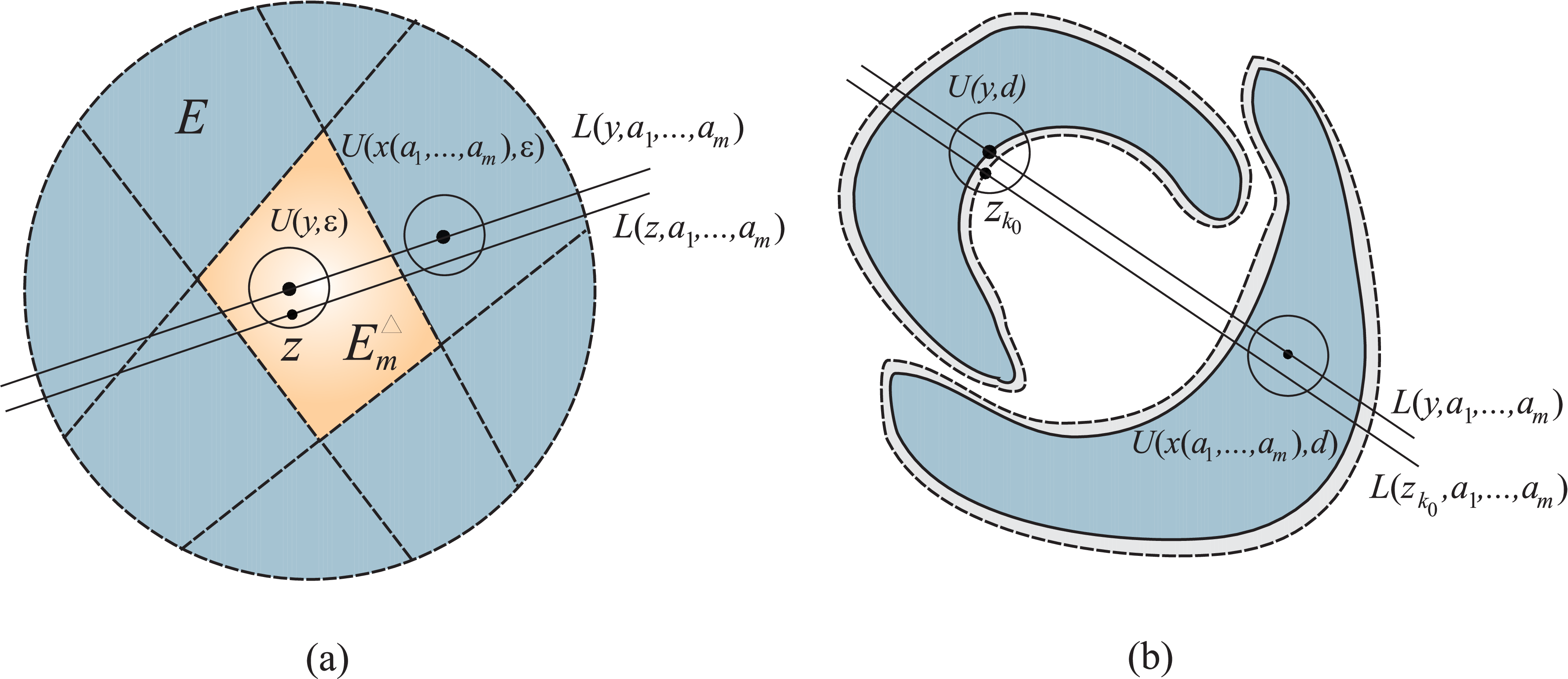}
\caption{}\label{Fig5}
\end{figure}

Let $\varepsilon:=\min\{\varepsilon_1,d\}$. Consider the neighborhood $U(y, \varepsilon)$ of the point $y$. Let $z\in U(y, \varepsilon)$, and $L(z,a_1,\ldots,a_m)$, $(a_1,\ldots,a_m) \in \mathcal L$, be an arbitrary  $m$-dimensional plane passing through $z$. Draw the $m$-dimensional plane $L(y,a_1,\ldots,a_m)$ parallel to $L(z,a_1,\ldots,a_m)$. Since $U(x(a_1,\ldots,a_m), \varepsilon)\subseteq U(x(a_1,\ldots,a_m), d)\subset E$ for the point $x(a_1,\ldots,a_m)$ correspondent to $ L (y,a_1,\ldots,a_m)$, it follows that $ L (z,a_1,\ldots,a_m)\cap U(x(a_1,\ldots,a_m), \varepsilon)\ne\varnothing$ and, therefore, $L (z,a_1,\ldots,a_m)\cap E\ne\varnothing$ for any $(a_1,\ldots,a_m) \in \mathcal L$; see Figure~\ref{Fig5}\,(a). Thus, $z$ is an $m$-nonconvexity point of $E$. Since $z$ is arbitrary, it implies that all points of $U(y, \varepsilon)$ are $m$-nonconvexity points of $E$. Hence, $y$ is an interior point of $E_m^{\triangle}$.
\end{proof}
\begin{corollary}
Let an open subset $E\subset \mathbb{R}^n$, $n\geqslant2$, belong to the class $\mathbf{WC^n_m}$, and $E_m^{\triangle}\ne\varnothing$, where $m=1,2,\ldots,n-1$. Then $E_m^{\triangle}$ is weakly $m$-convex.
\end{corollary}
\begin{proof} Since the set $E_m^{\triangle}$  is open, it follows that $\partial E_m^{\triangle}\cap E_m^{\triangle}=\varnothing$. Moreover, $\partial E_m^{\triangle}\cap E=\varnothing$, since $E$ is open and $E_m^{\triangle}\subset \mathbb{R}^n\setminus E$. Then, through every point $y\in\partial E_m^{\triangle}$, there passes an $m$-dimensional plane $L(y,a_1,\ldots,a_m)$, $(a_1,\ldots,a_m) \in \mathcal L$, that does not intersect $E$. In this case, $L(y,a_1,\ldots,a_m)\cap E_m^{\triangle}=\varnothing$ by Definition~\ref{def5}. Thus, $E_m^{\triangle}$ is weakly $m$-convex.
\end{proof}
\begin{theorem}
Let $E\subset \mathbb{R}^n$, $n\geqslant2$, be a closed subset such that $\mathrm{Int}\, E\ne\varnothing$. If  $E\in\mathbf{WC^n_m}$, then $\mathrm{Int}\, E\in \mathbf{WC^n_m}$ for every $m=1,2,\ldots,n-1$.
\end{theorem}
\begin{proof}
We prove the theorem by contradiction. Suppose that  $\mathrm{Int}\, E$ is not weakly $m$-convex. Then there exists a point $y\in\partial E\cap (\mathrm{Int}\,E)^{\triangle}_{m}$.

By Lemma~\ref{lem2}, for the point $y$ and every $(a_1,\ldots,a_m) \in \mathcal L$ there exist a point $x(a_1,\ldots,a_m) \in L (y,a_1,\ldots,a_m)\cap \mathrm{Int}\,E$ and a constant $d>0$ such that $U(x(a_1,\ldots,a_m),d)\subset \mathrm{Int}\,E$.

Now, consider the neighborhood $U(y, d)$ of the point $y$; see Figure~\ref{Fig5}\,(b).
Since $E$ is weakly $m$-convex, there is a family of open weakly $m$-convex sets $G_k$, $k=1,2,\ldots$, approximating $E$ from the outside. Then there exists an index $k_0\in \mathbb{N}$ such that $\partial G_{k_0}\cap U(y, d)\ne\varnothing$. Otherwise,  $U(y, d)\subset G_k$ for all $k\in\mathbb{N}$, since $y\in E\subset G_k$, and  $U(y, d)$ is connected and does not intersect $\partial G_{k}$. In this case, considering Definition~\ref{def2},  $U(y, d)\subset \bigcap\limits_kG_k=E$, which contradicts the fact that $y\in\partial E$.

Choose a point $z_{k_0}\in\partial G_{k_0}\cap U(y, d)$ and draw an arbitrary $m$-dimensional plane $L(z_{k_0},a_1,\ldots,a_m)$, $(a_1,\ldots,a_m) \in \mathcal L$,  through $z_{k_0}$. Consider the $m$-dimensional plane $L(y,a_1,\ldots,a_m)$ parallel to  $L(z_{k_0},a_1,\ldots,a_m)$.  Since $U(x(a_1,\ldots,a_m),d)\subset \mathrm{Int}\,E\subset E$ for the point $x(a_1,\ldots,a_m)$ corresponding to $L(y,a_1,\ldots,a_m)$, it follows that $L(z_{k_0},a_1,\ldots,a_m)\cap U(x(a_1,\ldots,a_m),d)\ne\varnothing$ and, therefore, $L(z_{k_0},a_1,\ldots,a_m)\cap E\ne\varnothing$. Since $G_{k_0}\supset E$, we also have $L(z_{k_0},a_1,\ldots,a_m)\cap G_{k_0}\ne\varnothing$.

 The $m$-dimensional plane $L(z_{k_0},a_1,\ldots,a_m)$ is arbitrary, so $z_{k_0}\in\partial G_{k_0}$ is an $m$-nonconvexity  point of $G_{k_0}$, which gives a contradiction with the fact that $G_{k_0}$ is weakly $m$-convex.
\end{proof}

\section*{Declarations}
This work was supported by the National Research Foundation of Ukraine in the framework of the research project 2025.07/0014 on the topic "Modern problems of Mathematical Analysis and Geometric Function Theory". The author declares no potential conflict of interest with respect to the research, authorship and publication of this article. All necessary data are included into the paper.


\begin{thebibliography}{9}

\bibitem{Ben}
Bendokat, T., Zimmermann, R. and Absil, P.-A. (2024).
A Grassmann manifold handbook: basic geometry and computational aspects.
\emph{Advances in Computational Mathematics}, \textbf{50}(1), Article 6.
https://doi.org/10.1007/s10444-023-10090-8.

\bibitem{Leh}
 Leichtwei{\ss}, K. (1980).
Konvexe Mengen,
VEB Deutscher Verlag der Wissenschaften, Berlin.

\bibitem{Osi2022_2}
Osipchuk, T.~M. (2022).
Topological and Geometric Properties of the Set of $1$-Nonconvexity
  Points of a Weakly $1$-Convex Set in the Plane.
\emph{Ukrainian Mathematical Journal}. \textbf{73}, 1918--1936.
https://doi.org/10.1007/s11253-022-02038-w.

\bibitem{Osi2023}
 Osipchuk, T.~M. and Tkachuk, M.~V. (2023).
On weakly 1-convex sets in the plane.
\emph{Proceedings of the International Geometry Center}. \textbf{16}(1), 42--49.
https://doi.org/10.15673/tmgc.v16i1.2440.

\bibitem{Osi2025}
Osipchuk, T.~M. (2025).
On weakly $1$-convex and weakly $1$-semiconvex sets in real Euclidean spaces.
\emph{Proceedings of the International Geometry Center}, \textbf{18}(1), 102--117.
https://doi.org/10.15673/pigc.v18i1.2947.

\bibitem{Osi2026_1}
Osipchuk, T.~M. (2026).
On weakly $(n-1)$-convex sets in low-dimensional Euclidean spaces.
\emph{Ukrains$\prime$kyi Matematychnyi Zhurnal}, \textbf{78}(7-8), 548--562.
https://doi.org/10.3842/umzh.v78i7-8.9669.

\bibitem{Osi2026_2}
Osipchuk, T.~M. (2026).
On weakly $2$-convex sets in three-dimensional Euclidean space.
\emph{Proceedings of the International Geometry Center}, \textbf{19}(2), paper 10, 11 pages.
https://doi.org/10.15673/pigc.3347.

\bibitem{Osi2026_3}
Osipchuk, T.~M. (2026).
On weakly $(n-1)$-convex sets and $(n-1)$-convex hulls in Euclidean spaces.
\emph{Matematychni Studii} (submitted for publication).

\bibitem{Zel01}
 Zelinskii, Yu.~B. and  Momot, I.~V. (2001).
On $(n, m)$-convex sets.
\emph{Ukrainian Mathematical Journal}. \textbf{53}, 482--487.
https://doi.org/10.1023/A:1012352607365.
\end{thebibliography}
\end{document}